\documentclass[12pt,a4paper,twoside]{article}
\usepackage{theorem}
\usepackage{amsfonts}
\usepackage{amsmath}
\usepackage{graphicx}

\newcommand{\R}{{\mathord{\mathbb R}}}

\newcommand{\N}{{\mathord{\mathbb N}}}
\newcommand{\C}{{\mathord{\mathbb C}}}

\newcommand{\mH}{{\mathcal H}}
\newcommand{\mB}{{\mathcal B}}

\newcommand{\mM}{{\mathcal M}}

\newcommand{\mN}{{\mathcal N}}

\newcommand{\f}{\varphi}
\newtheorem{thm}{Theorem}
\newtheorem{proposition}[thm]{Proposition}
\newtheorem{lemma}[thm]{Lemma}
\newtheorem{definition}[thm]{Definition}
{\theorembodyfont{\upshape} \newtheorem{rem}[thm]{\it Remark}}

\newtheorem{assumption}[thm]{Assumption}

\begin{document}
\pagestyle{myheadings}    
\markboth{W. H. Aschbacher}{Fully discrete Galerkin schemes for the nonlinear and nonlocal 
Hartree equation}
\title{Fully discrete Galerkin schemes for the nonlinear and nonlocal Hartree equation}
\author{Walter H. Aschbacher\footnote{aschbacher@ma.tum.de}\\ \\
Technische Universit\"at M\"unchen \\
Zentrum Mathematik, M5\\
85747 Garching, Germany
}
\date{}
\maketitle

\begin{abstract} 
We study the time dependent Hartree equation in the continuum, the 
semidiscrete, and the fully discrete setting. We prove 
existence-uniqueness, regularity, and approximation 
properties for the respective
schemes, and set the stage for a controlled numerical computation of 
delicate nonlinear and nonlocal features of the Hartree dynamics in 
various physical 
applications.
\end{abstract}



\section{Introduction}
In this paper, we study the nonlinear and nonlocal Hartree initial-boundary 
value problem for the (wave) function $\psi(x,t)$ being defined 
by\footnote{The notation ${\dot \psi}$ stands for the derivative of $\psi$ 
w.r.t. the time variable $t$.}
\begin{eqnarray}
\label{def:cp}
\left\{
\begin{array}{rl}
{\rm i}\,\dot\psi
=(-\Delta+v+\lambda V\ast |\psi|^2)\,\psi, & 
\,\,\mbox{if}\quad(x,t)\in \overline{\Omega}\times [0,T),\\
\psi=0, \hspace{4.0cm}     & 
\,\,\mbox{if}\quad (x,t)\in \partial\Omega\times [0,T),\\
\psi=\psi_0, \hspace{3.8cm}& 
\,\,\mbox{if}\quad (x,t)\in \overline{\Omega}\times\{0\},
\end{array}
\right.
\end{eqnarray}
where $\Omega$ is some domain in $\R^d$ with boundary 
$\partial\Omega$, and $T>0$ is some upper limit of the time interval on which 
we want to study the time evolution of $\psi$.
 Moreover, $v$ stands for an external potential, 
$\lambda$ denotes the coupling strength, and $V$ is the interaction potential 
responsible 
for the nonlinear and nonlocal interaction generated by the 
convolution.\footnote{To be made precise below.}
The system \eqref{def:cp} has many physical applications, in particular for 
the case 
$\Omega=\R^d$. As a first application, we mention the appearance of
\eqref{def:cp} within the context of the quantum mechanical description of
large systems of nonrelativistic bosons in their so-called mean field limit.
For the case of a local nonlinearity, i.e. $V=\delta$, an important 
application of equation \eqref{def:cp} lies in the domain  of Bose-Einstein 
condensation 
for repulsive interatomic forces where it governs the condensate wave 
function and is called the Gross-Pitaevskii equation. This dynamical 
equation, and its 
corresponding energy functional in the stationary case,
have been derived 
rigorously (see for example 
\cite{Spohn, Erdos} and \cite{Lieb}, respectively.\footnote{Moreover, the nonlocal 
Hartree equation has been derived for weakly
coupled fermions in \cite{Elgart}.}). In  \cite{Aschbacher},  
minimizers of this  nonlocal Hartree functional have been studied in the 
attractive case, and symmetry breaking has been established for sufficiently
large coupling. A large coupling phase segregation phenomenon has also been 
rigorously derived for a system of two coupled 
Hartree 
equations which are used to describe interacting Bose-Einstein condensates
(see \cite{Squassina, Aschbacher08} and references therein).
Such coupled systems also appear in the description of crossing sea states 
of weakly nonlinear dispersive surface water waves in hydrodynamics (see
for example \cite{Onorato,Shukla}), of  electromagnetic waves in a Kerr medium 
in nonlinear optics, and in nonlinear plasma physics. 
Furthermore, we would like to mention that
equation \eqref{def:cp} with attractive interaction potential 
possesses a so-called point particle
limit. Consider the situation where the initial condition is composed of 
several interacting 
Hartree minimizers sitting in an external potential which varies slowly on the 
length scale defined by the extension of the minimizers. It turns out that,
in a time regime inversely related to this scale, the center of mass of 
each minimizer follows a trajectory
which is governed, up to a small friction term, by Newton's equation of
motion for  interacting point particles in the slowly varying external 
potential. Hence,
in this limit, the system can be interpreted as the motion of interacting 
extended 
particles in a shallow external potential  and weakly coupled to  a dispersive
environment with which mass and energy can be exchanged through
 the friction term.
This allows to describe, and hence to numerically compute, some type of 
structure formation in Newtonian gravity (see \cite{Froehlich1, Froehlich2}).

The main content of the present paper consists in setting up the framework 
for the numerical analysis on bounded $\Omega$ which will be used in 
\cite{A09} for the study through
numerical computation of such phenomena, like, for example, the 
dissipation through
radiation for a Hartree minimizer oscillating in an external confining 
potential (see also \cite{diss}).

In Section \ref{sec:cp}, we start off with a brief study of the Hartree 
initial-value boundary problem \eqref{def:cp} in the  continuum setting
and we discuss its existence-uniqueness and regularity properties. In Section
\ref{sec:semidiscrete}, the system \eqref{def:cp} is discretized in space 
with the help of Galerkin theory. We derive existence-uniqueness and 
a bound on the $L^2$-approximation error. In the main Section \ref{sec:fd}, 
we proceed to the full discretization of
\eqref{def:cp}, more precisely, we discretize the foregoing semidiscrete 
problem in time focusing on two time discretization schemes of 
Crank-Nicholson type. The first is the so-called one-step one-stage 
Gauss-Legendre Runge-Kutta method which conserves the mass of the discretized
wave function under the discrete time evolution. The second one is 
the so-called Delfour-Fortin-Payre scheme which, besides the mass, also 
conserves the energy of the system. We prove existence-uniqueness 
using contraction methods  suitable for implementation in 
\cite{diss, A09}. Moreover, we derive a time quadratic accuracy estimate on the
$L^2$-error of these approximation schemes. In the proofs of these assertions,
we write down rather explicit expressions for the
bounds in order to have some qualitative
idea how to achieve a good numerical control of the fully discrete 
approximations of the Hartree initial-value boundary problem \eqref{def:cp}
for the computation of delicate nonlinear and nonlocal features of the various 
physical scenarios discussed above.

\section{The continuum problem}
\label{sec:cp}

As discussed in the Introduction, we  start off by briefly studying the 
Hartree initial-boundary value problem \eqref{def:cp} in a suitable 
continuum setting. For this purpose, we make the following assumptions
concerning the domain $\Omega$, the external potential $v$,
 and the interaction potential $V$,
a choice which is motivated by the perspective of the fully discrete 
problem and the numerical analysis dealt with
in Section \ref{sec:fd} and the numerical computations in 
\cite{A09}.\footnote{Some Hartree-dynamical computations have already been
performed in \cite{diss}.} 

\begin{assumption}
\label{ass:Omega}
$\Omega\subset\R^d$ is a bounded domain with smooth boundary 
$\partial\Omega$. 
\end{assumption}\vspace{-5mm}

\begin{assumption}
\label{ass:Vv}
Assumption \ref{ass:Omega} holds, $v\in C^\infty_0(\Omega,\R)$, and 
$V\in C^\infty_0(\R^d,\R)$.
\end{assumption}\vspace{-5mm}

\begin{assumption}
\label{ass:even}
Assumption \ref{ass:Vv} holds, and $V(-x)=V(x)$ for all $x\in\R^d$.
\end{assumption}

The Lebesgue and Sobolev spaces
used in 
the following are always defined over the domain $\Omega$ from Assumption
\ref{ass:Omega} unless something else is stated explicitly. Thus, we suppress
$\Omega$ in the notation of these spaces.  Moreover, under Assumption
\ref{ass:Vv}, let  the Hilbert space $\mH$,
the linear operator $A$  on $\mH$ with domain of definition $D(A)$, and the
nonlinear mapping $J$  on $\mH$ be given by
\begin{eqnarray*}
\mH&:=&L^2,\\
D(A)&:=&H^2\cap H^1_0,\\
A&:=&-\Delta,\\
J[\psi]&:=&v\psi+f[\psi],
\end{eqnarray*}
where the nonlinear mapping $f$ is defined by
\begin{eqnarray*}
f[\psi]&:=&\lambda g_V[|\psi|^2]\psi,\\ 
\lambda&\in& C^\infty_0(\Omega,\R),\\
g_V[\f](x)
&:=&\int_\Omega{\rm d}y\,\,V(x-y)\, \f(y).
\end{eqnarray*}

\begin{rem}
The function $\lambda$ stands for some space depending coupling function 
which can be 
chosen to be a smooth characteristic function of the domain $\Omega$. Such
a choice, on one hand, insures that all derivatives of $f[\psi]$ vanish at 
the boundary $\partial\Omega$, and, on the other hand, switches the nonlocal
interaction off in some neighborhood of $\partial\Omega$ where, in the 
numerical computation, transparent boundary conditions have to be matched with
the outgoing flow of $\psi$  (see \cite{diss, A09}).
\end{rem}

\begin{rem}
Under Assumption \ref{ass:Vv}, we have $g_V[\f]\in C^\infty_0(\R^d,\C)$ for
all $\f\in L^1$, and $f[\psi]\in L^2$ for all $\psi\in L^2$, see estimate
\eqref{est} below.
\end{rem}

We now make the following definition.

\begin{definition}
\label{def:cs}
Let Assumption \ref{ass:Vv} hold, and let $T\in (0,\infty)\cup \{\infty\}$.
We call a differentiable function $\psi:[0,T)\to L^2$ a continuum solution
of the Hartree initial-value problem \eqref{def:cp} with initial 
condition  $\psi_0\in D(A)$ if
\begin{eqnarray}
\label{cp}
\left\{
\begin{array}{l}
{\rm i}\dot\psi(t)
=A\psi(t)+J[\psi(t)],\quad
\forall t\in[0,T),\\
\,\psi(0)
=\psi_0.
\end{array}
\right. 
\end{eqnarray}
If $T<\infty$, the solution is called local, and if $T=\infty$, it is called
global.
\end{definition}

We make use of the following theorem to prove that there exists a 
unique global solution of the Hartree initial-value problem \eqref{def:cp} 
in the sense of Definition \ref{def:cs}.  In addition, this solution has 
higher regularity properties in time which are  required for the bounds on 
the constants appearing in the $L^2$-error 
estimates in the fully discrete setting of Section \ref{sec:fd}.

\begin{thm}(Cf. \cite[p.301]{RS2})
\label{thm:cs}
Let $\mH$ be a Hilbert space and $A$ a linear operator on $\mH$ with domain
of definition $D(A)$ and $A^\ast=A$. \\
(a) Let $n\in\N$, and let $J$ 
be a mapping which satisfies the following conditions for all 
$\psi,\xi\in D(A^k)$,
\begin{eqnarray}
\label{cond_a}
JD(A^k)
&\subseteq& D(A^k),\quad 
\forall k=1,...,n,\\
\label{cond_b}
\|J[\psi]\|
&\le& C(\|\psi\|) \|\psi\|,\\
\label{cond_c}
\|A^k J[\psi]\|
&\le& C(\|\psi\|,...,\|A^{k-1}\psi\|) \|A^k\psi\|,\quad 
\forall k=1,...,n,\\
\label{cond_d}
\|A^k(J[\psi]-J[\xi])\|
&\le& C(\|\psi\|,\|\xi\|,...,\|A^k\psi\|,\|A^k\xi\|)
\|A^k\psi-A^k\xi\|,\quad \forall k=0,...,n,
\end{eqnarray}
where each constant $C$ is a monotone increasing and everywhere finite function
of all its variables. Then, for each $\psi_0\in D(A^n)$, there exists a unique
local continuum solution in the sense of Definition \ref{def:cs} with
$\psi(t)\in D(A^n)$ for all $t\in[0,T)$.\\
(b) Moreover, this solution  is global if 
\begin{eqnarray}
\label{boundedness}
\|\psi(t)\|
\le C,\quad
\forall t\in[0,T).
\end{eqnarray}
(c) In addition to the conditions in (a), let $J$ satisfy the following 
conditions for all $1\le k< n$: if
$\psi\in C^k([0,\infty),\mH)$ with $\frac{{\rm d}^l}{{\rm d} t^l}\psi(t) 
\in D(A^{n-l})$ for all $0\le l\le k$ and   
$A^{n-l}\frac{{\rm d}^l}{{\rm d} t^l}\psi(t)\in C([0,\infty),\mH)$ 
for all $0\le l\le k$, then 
\begin{eqnarray}
\label{cond_e}
J[\psi(t)]
&& \mbox{\hspace{-7mm}is $k$ times differentiable},\\
\label{cond_f}
\frac{{\rm d}^k}{{\rm d} t^k}J[\psi(t)]
\!\!&\in&\!\! D(A^{n-k-1}),\\
\label{cond_g}
A^{n-k-1}\frac{{\rm d}^k}{{\rm d} t^k}J[\psi(t)]
\!\!&\in&\!\! C([0,\infty),\mH).
\end{eqnarray}
If this condition holds, then the local solution from (a) is $n$ times 
differentiable  and  $\frac{{\rm d}^k}{{\rm d} t^k}\psi(t) \in D(A^{n-k})$ 
for all $1\le k\le n$.
\end{thm}

\begin{rem}
With the help  of Theorem \ref{thm:cs} (and its proof) the constants in 
the estimates \eqref{cfd:c4} and \eqref{ifd:c4} below 
on the $L^2$-error of 
the fully discrete approximations are finite and can be estimated 
explicitly.
\end{rem}

{\it Proof}\, Let us start off by checking condition \eqref{cond_b}.

{\it Condition} \eqref{cond_b}:\quad ${\|J[\psi]\|}_{L^2}\le C({\|\psi\|}_{L^2}) 
{\|\psi\|}_{L^2}$\\
Using Assumption \ref{ass:Vv} and the estimate 
\begin{eqnarray}
\label{est}
{\|\lambda g_V[\f\psi]\chi\|}_{L^2}
\le {\|\lambda\|}_{L^\infty} {\|V\|}_{L^\infty(\R^d)} 
{\|\f\|}_{L^2}{\|\psi\|}_{L^2}{\|\chi\|}_{L^2},
\end{eqnarray}
we immediately get
\begin{eqnarray*}
{\|J(\psi)\|}_{L^2}
\le ({\|v\|}_{L^\infty}
+{\|\lambda\|}_{L^\infty} {\|V\|}_{L^\infty(\R^d)}{\|\psi\|}_{L^2}^2)
{\|\psi\|}_{L^2},
\end{eqnarray*}
and the prefactor $C(\|\psi\|)$ from \eqref{cond_b} is a monotone increasing 
function of ${\|\psi\|}_{L^2}$. Next, let us check condition \eqref{cond_c}.

{\it Condition} \eqref{cond_c}:\quad ${\|\Delta^k J[\psi]\|}_{L^2} 
\le C({\|\psi\|}_{L^2}){\|\Delta^k\psi\|}_{L^2}$ for all $k\in\N$\\
In order to show \eqref{cond_c}, we have to control the $L^2$-norm of 
$\Delta^k(v\psi)$
and of $\Delta^k(\lambda g_V[|\psi|^2]\psi)$. Hence, we write the powers of 
the Dirichlet-Laplacian as follows,
\begin{eqnarray}
\label{Deltak1}
\Delta^k(\f\psi)
&=&\hspace{-1cm}\sum_{\substack{j_1,...,j_k=1,...,d\\
\alpha_1,\beta_1,...,\alpha_k,\beta_k=0,1}}
\hspace{-5mm}c_{\alpha_1\beta_1...\alpha_k\beta_k}\,
(\partial_{j_1}^{\alpha_1}\partial_{j_1}^{\beta_1}\!...
\partial_{j_k}^{\alpha_k}\partial_{j_k}^{\beta_k}\f)
(\partial_{j_1}^{1-\alpha_1}\partial_{j_1}^{1-\beta_1}\!\!...
\partial_{j_k}^{1-\alpha_k}\partial_{j_k}^{1-\beta_k}\psi),\\
\label{Deltak2}
\Delta^k(g_V[\f]\psi)
&=&\hspace{-1cm}\sum_{\substack{j_1,...,j_k=1,...,d\\
\alpha_1,\beta_1,...,\alpha_k,\beta_k=0,1}}
\hspace{-5mm}c_{\alpha_1\beta_1...\alpha_k\beta_k}\,
g_{\partial_{j_1}^{\alpha_1}\partial_{j_1}^{\beta_1}\!...
\partial_{j_k}^{\alpha_k}\partial_{j_k}^{\beta_k}V}[\f]\,
(\partial_{j_1}^{1-\alpha_1}\partial_{j_1}^{1-\beta_1}\!\!...
\partial_{j_k}^{1-\alpha_k}\partial_{j_k}^{1-\beta_k}\psi),
\end{eqnarray}
where $c_{\alpha_1\beta_1...\alpha_k\beta_k}$ denote some combinatorial
constants. Hence, using \eqref{est} and the following Schauder type 
estimate,\footnote{I.e. the elliptic regularity estimate of generalized 
solutions
up to the boundary, see for example \cite[p.383]{Zeidler}.}
\begin{eqnarray*}
{\|\psi\|}_{H^{2+m}}
\le C {\|\Delta\psi\|}_{H^m},\quad
\forall \psi\in H^{2+m}\cap H^1_0,\quad
\forall m\in\N_0,
\end{eqnarray*}
we get the following bounds on \eqref{Deltak1} and \eqref{Deltak2},
\begin{eqnarray}
\label{est1}
{\|\Delta^k(\f\psi)\|}_{L^2}
&\le& C {\|\f\|}_{W^{2k,\infty}} {\|\psi\|}_{H^{2k}}\nonumber\\
&\le& C {\|\f\|}_{W^{2k,\infty}} {\|\Delta^k\psi\|}_{L^2},\\
\label{est2}
{\|\Delta^k(g_V[\f]\psi)\|}_{L^2}
&\le& C {\|V\|}_{W^{2k,\infty}(\R^d)} {\|\f\|}_{L^1}
{\|\Delta^k\psi\|}_{L^2}.
\end{eqnarray}
Using \eqref{est1}, and \eqref{est2} twice, we arrive at 
\begin{eqnarray}
\label{cond_cc}
{\|\Delta^k J[\psi]\|}_{L^2} 
\le C ({\|v\|}_{W^{2k,\infty}}
+{\|\lambda\|}_{W^{2k,\infty}} {\|V\|}_{W^{2k,\infty}(\R^d)} 
{\|\psi\|}_{L^2}^2) {\|\Delta^k\psi\|}_{L^2}, 
\end{eqnarray}
where the prefactor $C(\|\psi\|, \|A\psi\|,...,\|A^{k-1}\psi\|)$ from
\eqref{cond_c} depends on ${\|\psi\|}_{L^2}$ only and is monotone increasing 
in ${\|\psi\|}_{L^2}$. Next, we check condition \eqref{cond_a}.

{\it Condition} \eqref{cond_a}:\quad $JD(\Delta^k)\subseteq D(\Delta^k)$ 
for all 
$k\in\N$\\
Due to \eqref{cond_cc} and since  $D(\Delta^k)=\{\psi\in H^2\cap H^1_0\,|\,
\Delta\psi\in H^2\cap H^1_0,...,\Delta^{k-1}\psi\in H^2\cap H^1_0\}$, 
it remains to show that $\Delta^{j}J(\psi)\in H^1_0$ for all 
$\psi\in D(\Delta^k)$ and for all $j=0,...,k$. But this follows since
$v,\lambda\in C^\infty_0$ and from the fact that $C^\infty(\overline{\Omega})$ 
is dense in $H^m$ w.r.t. the $H^m$-norm for all $m\in\N_0$.
Let us next check condition \eqref{cond_d}.

{\it Condition} \eqref{cond_d}:\quad 
${\|\Delta^k(J[\psi]-J[\xi])\|}_{L^2}
\le C({\|\psi\|}_{L^2},{\|\xi\|}_{L^2}, {\|\Delta^k\psi\|}_{L^2},
{\|\Delta^k\xi\|}_{L^2}) {\|\Delta^k\psi-\Delta^k\xi\|}_{L^2}$\\
In order to show \eqref{cond_d}, we write the difference with the help
of the decomposition
\begin{eqnarray}
\label{dc1}
g_V[\f_1\psi_1]\chi_1
-g_V[\f_2\psi_2]\chi_2
&=&g_V[\f_1(\psi_1-\psi_2)]\chi_1\nonumber\\
&&+\,g_V[(\f_1-\f_2)\psi_2]\chi_1\nonumber\\
&&+\,g_V[\f_2\psi_2](\chi_1-\chi_2).
\end{eqnarray}
Each term on the r.h.s. of \eqref{dc1} can then be estimated with the help
of \eqref{est2}. Hence, as in \eqref{cond_c}, we get
\begin{eqnarray}
\label{cond_dd}
{\|\Delta^k(J[\psi]-J[\xi])\|}_{L^2}
&\le& C ({\|v\|}_{W^{2k,\infty}}
+{\|\lambda\|}_{W^{2k,\infty}} {\|V\|}_{W^{2k,\infty}(\R^d)}\nonumber\\
&&\hspace{7mm}\cdot ({\|\psi\|}_{L^2}^2
+({\|\psi\|}_{L^2}+{\|\xi\|}_{L^2}){\|\Delta^k\psi\|}_{L^2}))\,
{\|\Delta^k\psi-\Delta^k\xi\|}_{L^2},
\end{eqnarray}
where the prefactor $C(\|\psi\|,\|\xi\|,...,\|A^k\psi\|,\|A^k\xi\|)$ 
from \eqref{cond_d} depends on the lowest and the highest power of $\Delta$
only, and it is monotone increasing in all its variables.

{\it Condition} \eqref{boundedness}:\quad 
${\|\psi(t)\|}_{L^2}={\|\psi_0\|}_{L^2}$ for all $t\in [0,T)$\\
This condition is satisfied due to Proposition \ref{prop:conservation} below.

{\it Condition} \eqref{cond_e}:\quad $J[\psi(t)]$ is $k$ times 
differentiable in $L^2$ w.r.t. time $t$\\
Let $n\in\N$ be fixed, and let $k=1$. Then, it is shown in \cite[p.299]{RS2}
that part (a) and conditions \eqref{cond_e}, \eqref{cond_f}, 
and \eqref{cond_g} for $k=1$ imply that
$\psi\in C^2([0,\infty),L^2)$ with $\frac{{\rm d}^2}{{\rm d} t^2}\psi(t) 
\in D(A^{n-2})$ and   
$A^{n-2}\frac{{\rm d}^2}{{\rm d} t^2}\psi(t)\in C([0,\infty),L^2)$. Then, using
 conditions \eqref{cond_e}, \eqref{cond_f}, 
and \eqref{cond_g} for  subsequent $1\le k<n$
leads to the  claim of part (c) by iteration. Hence, we have to verify that
conditions \eqref{cond_e}, \eqref{cond_f}, 
and \eqref{cond_g} are satisfied for $1\le k<n$. To this end, we make use of 
decomposition \eqref{dc1} to exemplify the case $k=1$ and to note that the
cases for $k\ge 2$ are analogous. In order to show that $J[\psi(t)]$ is
differentiable in $L^2$, we write, using \eqref{dc1},
\begin{eqnarray*}
\frac{J[\psi(t+h)]-J[\psi(t)]}{h}
&=&v\,\frac{\psi(t+h)-\psi(t)}{h}\\
&&+\lambda g_V[\psi(t+h)\tfrac{\bar\psi(t+h)-\bar\psi(t)}{h}]\psi(t+h)\\
&&+\lambda g_V[\tfrac{\psi(t+h)-\psi(t)}{h}\,\bar\psi(t+h)]\psi(t+h)\\
&&+\lambda g_V[|\psi(t)|^2] \frac{\psi(t+h)-\psi(t)}{h}.
\end{eqnarray*}
Applying \eqref{dc1} and \eqref{est}, we find that 
$J[\psi(t)]$ is 
differentiable in $L^2$ w.r.t. $t$ with derivative
\begin{eqnarray}
\label{deriv}
\frac{{\rm d}}{{\rm d}t}J[\psi(t)]
=v\dot\psi(t)
+\lambda g_V[\psi(t){\dot{\bar\psi}}(t)]\psi(t)
+\lambda g_V[\dot\psi(t)\bar\psi(t)]\psi(t)
+\lambda g_V[|\psi(t)|^2] \dot\psi(t).
\end{eqnarray}
Making use of \eqref{dc1}, \eqref{est1}, and \eqref{est2} in the estimate
of ${\|\Delta^{n-2}(\tfrac{{\rm d}}{{\rm d}t}J[\psi(t)]
-\tfrac{{\rm d}}{{\rm d}t}J[\psi(s)])\|}_{L^2}$ similarly to 
\eqref{cond_dd}, we find that 
$\tfrac{{\rm d}}{{\rm d}t}J[\psi(t)]\in D(\Delta^{n-2})$ and that 
$\Delta^{n-2}\tfrac{{\rm d}}{{\rm d}t}J[\psi(t)]\in C([0,\infty),L^2)$. Due to
the structure of \eqref{deriv}, we can iterate the foregoing procedure to 
arrive at the assertion.
\hfill$\Box$\vspace*{5mm}

In order to verify condition \eqref{boundedness}, we define 
the  mass $\mM[\psi]$ and energy $\mH[\psi]$ of a function $\psi\in H^1$ by
\begin{eqnarray*}
\mM[\psi]
&:=&{\|\psi\|}_{L^2}^2,\\
\mH[\psi]
&:=&{\|\nabla\psi\|}_{L^2}^2
+{(\psi,v\psi)}_{L^2}
+\frac12\,{(\psi,f[\psi])}_{L^2}.
\end{eqnarray*}

We then have the following.

\begin{proposition}
\label{prop:conservation}
Let Assumption \ref{ass:even} hold, and let $\psi$ be the unique local 
continuum solution of  Theorem \ref{thm:cs}. 
Then,  the  mass 
and the energy of $\psi$ are conserved under the time evolution, 
\begin{eqnarray}
\label{cp:conservation}
\mM[\psi(t)]
&=&\mM[\psi_0],\quad 
\forall t\in [0,T),\\
\mH[\psi(t)]
&=&\mH[\psi_0],\quad 
\forall t\in [0,T).\nonumber
\end{eqnarray}
\end{proposition}

{\it Proof}\, From  Theorem \ref{thm:cs} it follows that
the
function $t\mapsto \mM[\psi(t)]$ belongs to $C^1([0,T),\R_0^+)$ with
\begin{eqnarray*}
\frac{{\rm d}}{{\rm d}t}\,\mM[\psi(t)]
=2\,{\rm Re}\, {(\dot \psi(t),\psi(t))}_{L^2},
\end{eqnarray*}
which vanishes due to \eqref{cp}. For the conservation of the 
energy, we have the following three parts. First, using the regularity of 
$\psi$ in time and 
$(\psi,-\Delta \psi)_{L^2}={\|\nabla\psi\|}_{L^2}^2$ for all 
$\psi\in H^2\cap H^1_0$, we observe that  
$t\mapsto {\|\nabla\psi(t)\|}_{L^2}^2$
 belongs to $C^1([0,T),\R_0^+)$ and has the derivative
\begin{eqnarray}
\label{deriv1}
\frac{{\rm d}}{{\rm d}t}\,{\|\nabla\psi(t)\|}_{L^2}^2
=2\,{\rm Re}\,{(\dot\psi(t),-\Delta\psi(t))}_{L^2}.
\end{eqnarray}
Second, the function $t\mapsto {(\psi(t),v\psi(t))}_{L^2}$ 
belongs to $C^1([0,T),\R)$ and has the derivative 
\begin{eqnarray}
\label{deriv2}
\frac{{\rm d}}{{\rm d}t}\,{(\psi(t),v\psi(t))}_{L^2}
=2\,{\rm Re}\,{(\dot \psi(t),v\psi(t))}_{L^2}.
\end{eqnarray}
Third, using $|\psi|^2-|\f|^2=\psi(\bar\psi-\bar\f)+(\psi-\f)\bar\f$ in the
decomposition of $f[\psi]-f[\f]$ as in  \eqref{dc1}, we get 
$\frac{{\rm d}}{{\rm d}t}g_V[|\psi|^2]\psi
=g_V[|\psi|^2]\dot\psi
+2\,{\rm Re}\,(g_V[\bar\psi\dot\psi])\,\psi$
in $L^2$, and therefore the function $t\to {(\psi(t),f[\psi(t)])}_{L^2}$ 
belongs to $C^1([0,T),\R)$ and has derivative
\begin{eqnarray}
\label{deriv3}
\frac{{\rm d}}{{\rm d}t}\,{(\psi(t),f[\psi(t)])}_{L^2}
=4\,{\rm Re}\,{(\dot \psi(t),f[\psi(t)])}_{L^2},
\end{eqnarray}
where we used Assumption \ref{ass:even} to write 
${(g[\bar\psi\dot\psi]\psi,\psi)}_{L^2}
={(\dot\psi,f[\psi])}_{L^2}$. Finally, if we take the 
scalar product of \eqref{cp} with $\dot\psi$ and the real part of the
resulting equation, we get
\begin{eqnarray}
\label{ScalarProduct}
0={\rm Re}\,\,{\rm i}{\|\dot\psi(t)\|}_{L^2}^2
={\rm Re}\,\big[{(\dot\psi(t),-\Delta\psi(t))}_{L^2}
+{(\dot\psi(t),v\psi(t))}_{L^2}
+{(\dot\psi(t),f[\psi(t)])}_{L^2}\big].
\end{eqnarray}
Plugging  
\eqref{deriv1}, \eqref{deriv2}, and \eqref{deriv3} into \eqref{ScalarProduct},
we find the conservation of the energy $\mH[\psi(t)]$.
\hfill$\Box$\vspace{5mm}

\begin{rem}
For more general interaction potentials 
$V$, in particular in the local
case $f[\psi]=|\psi|^2\psi$ in $d=2$,\footnote{We are mainly interested in
$d=2$ for the numerical computations.} one can use an estimate from 
\cite{Brezis} which 
controls the 
$L^\infty$-norm of a function $\psi\in H^1$ by the square root of the 
logarithmic growth of the $H^2$-norm,
\begin{eqnarray*}
{\|\psi\|}_{L^\infty}
\le C\,\big(1+\sqrt{\log(1+{\|\psi\|}_{H^2})}\big),
\end{eqnarray*}
where the constant $C$ depends on ${\|\psi\|}_{H^1}$.
This estimate allows to bound the graph norm of the continuum
solution by a double
exponential growth, and, hence, makes the solution global.
\end{rem}

Taking the $L^2$-scalar product of \eqref{cp} w.r.t. functions 
$\f\in H^1_0$, and using again that 
${(\f,-\Delta\psi)}_{L^2}={(\nabla\f,\nabla\psi)}_{L^2}$ for all 
$\psi\in H^2\cap H^1_0$ and for all $\f\in H^1_0$, we get the following
weak formulation of the continuum problem \eqref{cp},
\begin{eqnarray}
\label{def:wcp}
\left\{
\begin{array}{ll}
{\rm i} {(\f,{\dot \psi})}_{L^2}
={(\nabla\f,\nabla\psi)}_{L^2}
+{(\f,v\psi)}_{L^2}
+{(\f,f[\psi])}_{L^2},\quad
\forall \f\in H^1_0,\, \,\forall t\in[0,T),\\
\hspace{7.5mm}\psi(0)=\psi_0.
\end{array}
\right.
\end{eqnarray}

The formulation \eqref{def:wcp} is the starting point for 
a suitable discretization in space of the original continuum problem.
We will discuss such a semidiscrete approximation in  Section 
\ref{sec:semidiscrete}.

\section{The semidiscrete approximation}
\label{sec:semidiscrete}

In this section, we discretize the problem \eqref{def:wcp} in space
with the help of  Galerkin theory which makes use of a family 
$\{S_h\}_{h\in (0,1)}$
of finite dimensional subspaces approximating the infinite dimensional
problem in the following precise sense.

\begin{assumption}
\label{ass:Sh}
The family $\{S_h\}_{h\in (0,1)}$ of subspaces of $H_0^1$ 
has the property
\begin{eqnarray*}
S_h\subset C(\overline{\Omega})\cap H_0^1, 
\quad \dim S_h=N_h<\infty,
\quad \forall h\in (0,1).
\end{eqnarray*}
\end{assumption}

\begin{rem}
\label{rem:Lagrange}
For the numerical computation in \cite{A09}, the physical space 
is (a smoothly bounded superset of) the open square 
$\Omega=(0,D)^2\subset\R^2$ with $D>0$ whose closure is 
the union
of the $(n-1)^2$ congruent closed subsquares generated by dividing each side 
of $\Omega$ equidistantly into $n-1$ intervals. Let us  
denote by $N_h=(n-2)^2$ the total number of interior vertices of  
this  
lattice and by $h=D/(n-1)$ the lattice spacing.\footnote{\label{foot:bijection}As  
bijection from the one-dimensional to the two-dimensional  
lattice numbering, we may use the mapping  
$\tau:\{0,...,m-1\}^{\times 2} 
\to \{0,...,m^2-1\}$ with $j=\tau(m_1,m_2):=m_1+m_2m$.}  
Moreover, let us   
choose the  Galerkin space $S_h$ to be spanned by the bilinear  
Lagrange rectangle finite elements  
$\varphi_j\in C(\overline{\Omega})$ whose reference  
basis function $\varphi_0:\overline{\Omega}\to [0,\infty)$  
is defined on  
its support $[0,2h]^{\times 2}$ by 
\begin{eqnarray} 
\label{def:blfe} 
\varphi_0(x,y) 
:=\frac{1}{h^2}\left\{ 
\begin{array}{rl} 
xy,           & \mbox{if } (x,y)\in [0,h]^{\times 2},\\ 
(2h-x)y,      & \mbox{if } (x,y)\in [h,2h]\times [0,h],\\ 
(2h-x)(2h-y), & \mbox{if } (x,y)\in [h,2h]^{\times 2},\\ 
x(2h-y),      & \mbox{if } (x,y)\in [0,h]\times [h,2h], 
\end{array} 
\right. 
\end{eqnarray} 
see Figure \ref{fig:refel}. The functions $\varphi_j$ are then defined to 
be of the form \eqref{def:blfe} having their support translated by 
$(m_1 h,m_2 h)$ with $m_1,m_2=0,...,m-3$. Hence, with this choice, 
we have $S_h\subset C(\overline{\Omega})\cap H_0^1(\Omega)$ and 
$\dim S_h=N_h$.

\begin{figure}[h!] 
\begin{center} 
\includegraphics[width=4cm,height=4.5cm]{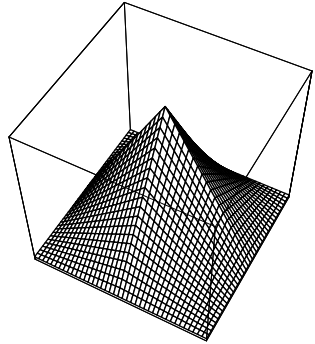}
\caption{$\varphi_0(x,y)$ on its support $[0,2h]^{\times 2}$ with  
maximum at vertex $(h,h)$.} 
\label{fig:refel} 
\end{center} 
\end{figure}

\end{rem}

Motivated by the weak formulation \eqref{def:wcp}, we make the following 
definition.

\begin{definition}
\label{def:sdp}
Let Assumptions \ref{ass:Vv} and \ref{ass:Sh} hold. 
We call $\psi_h: [0,T)\to S_h$ with $\psi_h,\dot\psi_h\in L^2(0,T;S_h)$ 
a semidiscrete solution of the Hartree initial 
boundary-value problem \eqref{def:cp}  with initial condition 
$\psi_{0h}\in S_h$ if
\begin{eqnarray}
\label{sdp}
\left\{
\begin{array}{l}
{\rm i}\, \frac{{\rm d}}{{\rm d}t}\,{(\f,\psi_h)}_{L^2}
={(\nabla\f,\nabla\psi_h)}_{L^2}
+{(\f,v\psi_h)}_{L^2}
+{(\f,f[\psi_h])}_{L^2},\quad 
\forall\f\in S_h,\,\, \forall t\in [0,T),\\ 
\hspace{12.2mm}\psi_h(0)=\psi_{0h}.
\end{array}
\right.
\end{eqnarray}
\end{definition}\vspace{0mm}

\begin{rem}
\label{rem:Gelfand}
In general, the weak problem \eqref{def:wcp} is set up using the 
Gelfand evolution triple 
$H_0^1\subset L^2\subset {(H_0^1)}^\ast=H^{-1}$. One then looks for weak 
solutions $\psi\in W_2^1(0,T; H_0^1, L^2)\subset C([0,T),L^2)$ 
motivating Definition \ref{def:sdp}.
\end{rem}

We assume the Galerkin subspace $S_h$ from Assumption \ref{ass:Sh} to 
satisfy the following additional approximation and inverse inequalities. 

\begin{assumption}
\label{ass:AccuracyOrder}
Let Assumption \ref{ass:Sh} hold. Then, there exists 
a constant $C_A>0$ s.t. 
\begin{eqnarray*}
\inf_{\f\in S_h}\big({\|\psi-\f\|}_{L^2}
+h{\|\psi-\f\|}_{H^1}\big)
\le C_A h^2 {\|\psi\|}_{H^2},
\quad\forall \psi\in H^2\cap H^1_0.
\end{eqnarray*}
\end{assumption}

\begin{rem}
For an order of accuracy $r\ge 2$ of the family $\{S_h\}_{h\in (0,1)}$,
the usual assumption replaces the r.h.s. by $C_A h^s {\|\psi\|}_{H^s}$
and is asked to hold for all $\psi\in H^s\cap H^1_0$. 
For simplicity, 
we stick to Assumption \ref{ass:AccuracyOrder}.
\end{rem}

\begin{assumption}
\label{ass:inverse}
Let Assumption \ref{ass:Sh} hold. Then, there exists a constant $C_B>0$ s.t. 
\begin{eqnarray*}
{\|\f\|}_{H^1}\le C_B h^{-1} {\|\f\|}_{L^2},
\quad\forall \f\in S_h.
\end{eqnarray*}
\end{assumption}

\begin{rem}
For the two-dimensional bilinear Lagrange finite element setting of Remark 
\ref{rem:Lagrange}, both Assumption 
\ref{ass:AccuracyOrder} and Assumption \ref{ass:inverse} 
hold.\footnote{See for example  \cite[p.109,111]{Brenner}.}
\end{rem}

Furthermore, we make an assumption on the approximation quality of the
initial condition $\psi_{0h}\in S_h$ of the semidiscrete problem 
\eqref{sdp}
compared to the initial condition $\psi_0\in H^2\cap H^1_0$ of the continuum 
problem \eqref{cp}.
\begin{assumption}
\label{ass:initial}
Let Assumption \ref{ass:Sh} hold. Then, there exists a constant $C_0>0$ s.t.
\begin{eqnarray}
\label{def:initial}
{\|\psi_0-\psi_{0h}\|}_{L^2}
\le C_0 h^2.
\end{eqnarray}
\end{assumption}

The semidiscrete scheme has the following conservation properties.

\begin{proposition}
\label{sdp:ExUnique}
Let Assumptions \ref{ass:even} and \ref{ass:Sh} hold, and let
$\psi_h$ be a semidiscrete solution of the Hartree initial boundary-value 
problem \eqref{def:cp} in the sense of Definition \ref{def:sdp}. 
Then, the mass and energy of $\psi_h$ are conserved under the time evolution,
\begin{eqnarray}
\label{sdp:conservation}
\mM[\psi_h(t)]=\mM[\psi_{0h}],\quad
\forall t\in[0,T),\\
\mH[\psi_h(t)]=\mH[\psi_{0h}],\quad
\forall t\in[0,T).\nonumber
\end{eqnarray}
\end{proposition}

{\it Proof}\, If we plug $\f=\psi_h(t)$ into 
\eqref{sdp} and take the imaginary part of the resulting equation,
we get the conservation of the mass. If
we plug $\f=\dot\psi_h(t)$ into 
\eqref{sdp} and take the real part of the resulting equation, we get
the conservation of the energy using Assumption \ref{ass:even}.\footnote{In
the sense of Remark \ref{rem:Gelfand}.}
\hfill$\Box$\vspace{5mm} 

Existence-uniqueness is addressed in the following.

\begin{thm}
\label{sdp:ExUniqueCons}
Let Assumptions \ref{ass:even} and \ref{ass:Sh} hold. Then, there exists a
unique global semidiscrete solution $\psi_h$ of the Hartree 
initial-boundary value problem \eqref{def:cp} in the sense of Definition 
\ref{def:sdp}. 
\end{thm}

{\it Proof}\, Let $\{\f_j\}_{j=1}^{N_h}$ be a basis of the Galerkin space 
$S_h$, and let us write
\begin{eqnarray}
\label{expansion1}
\psi_h(t)
=\sum_{j=1}^{N_h}z_j(t)\,\f_j.
\end{eqnarray}
Plugging \eqref{expansion1} into the semidiscrete system \eqref{sdp},
we get for $z(t):=(z_1(t),...,z_{N_h}(t))\in\C^{N_h}$, 
\begin{eqnarray}
\label{sdp2}
\left\{
\begin{array}{l}
{\rm i}\,\dot z(t)
=A^{-1}(B+Y)z(t)
+A^{-1}H[z(t)]z(t),
\quad\forall t\in[0,T),\\
z(0)=z_0,
\end{array}
\right.
\end{eqnarray}
where $\psi_{0h}=\sum_{j=1}^{N_h}(z_0)_j\,\f_j$ and the matrices 
$A,B\in\C^{N_h\times N_h}$ are the positive definite mass
and stiffness matrices, respectively,
\begin{eqnarray*}
A_{ij}&:=&{(\f_i,\f_j)}_{L^2},\\
B_{ij}&:=&{(\nabla\f_i,\nabla\f_j)}_{L^2}.
\end{eqnarray*}
Moreover, $Y\in\C^{N_h\times N_h}$ is the external potential matrix,
\begin{eqnarray*}
Y_{ij}:={(\f_i,v\f_j)}_{L^2}, 
\end{eqnarray*}
and the matrix-valued function $H:\C^{N_h}\to\C^{N_h\times N_h}$ is defined by 
\begin{eqnarray*}
{H[z]}_{ij}
:=\sum_{k,l=1}^{N_h}\bar z_k z_l\, {(\f_i,\lambda g[\bar\f_k\f_l]\f_j)}_{L^2}.
\end{eqnarray*}
Since the function $\C^{N_h}\ni z\mapsto A^{-1}(B+Y)z+A^{-1}H[z]z\in\C^{N_h}$
is locally Lipschitz continuous analogously to the continuum case, 
the Picard-Lindel\"of theory for ordinary differential equations implies
local existence and uniqueness of the initial value problem \eqref{sdp2}.
Moreover, this local solution is a global solution if it remains restricted 
to a compact subset of $\C^{N_h}$. But this is the case due to the mass 
conservation from
\eqref{sdp:conservation}. 
\hfill$\Box$\vspace{5mm} 

We next turn to the $L^2$-error estimate of the semidiscretization. For that
purpose, we introduce the Ritz projection.\footnote{Also called elliptic 
projection.}

\begin{definition}
\label{def:Ritz}
Let Assumption \ref{ass:Sh} hold. Then, 
the Ritz projection $R_h: H^1_0\to S_h$ is defined to be the 
orthogonal projection from $H^1_0$ onto $S_h$ w.r.t. the 
Dirichlet scalar product ${(\nabla\cdot,\nabla\cdot)}_{L^2}$ on 
$H^1_0$,  i.e.
\begin{eqnarray}
\label{Ritz}
{(\nabla\f,\nabla R_h\psi)}_{L^2}
={(\nabla\f,\nabla\psi)}_{L^2},
\quad\forall \f\in S_h.
\end{eqnarray}
\end{definition}

The Ritz projection satisfies the following error estimate.

\begin{lemma}{(Cf. \cite[p.8]{VT97})}
Let Assumptions \ref{ass:Sh} and \ref{ass:AccuracyOrder} hold. 
Then, there exists a constant $C_R>0$ s.t.
\begin{eqnarray}
\label{RitzEstimate}
{\|(1-R_h)\psi\|}_{L^2}+h{\|(1-R_h)\psi\|}_{H^1}
\le C_R h^2 {\|\psi\|}_{H^2},
\quad\forall \psi\in H^2\cap H^1_0.
\end{eqnarray}
\end{lemma}

The next theorem is the main assertion of this section. 

\begin{thm} 
\label{thm:ErrorSD}
Let  Assumptions \ref{ass:Vv}, \ref{ass:Sh},
\ref{ass:AccuracyOrder}, and \ref{ass:initial}  hold,
and let $\psi$ be the solution of the continuum problem from Theorem 
\ref{thm:cs} and $\psi_h$ the solution of the semidiscrete
problem from Theorem \ref{sdp:ExUniqueCons}. Then, for any $0<T<\infty$,
there exists a constant $C_E>0$ s.t. 
\begin{eqnarray*}
\max_{t\in[0,T]}{\|\psi(t)-\psi_h(t)\|}_{L^2}
\le C_E h^2.
\end{eqnarray*}
\end{thm}

{\it Proof}\, We decompose the difference of $\psi$ and $\psi_h$ as
\begin{eqnarray*}
\psi(t)-\psi_h(t)=\rho(t)+\theta(t),
\end{eqnarray*}
where $\rho(t)$ and $\theta(t)$ are defined with the 
help of the Ritz projection $R_h$ from \eqref{Ritz} by
\begin{eqnarray*}
\rho(t):=(1-R_h)\psi(t),\quad 
\theta(t):=R_h\psi(t)-\psi_h(t).
\end{eqnarray*}
Making use of the schemes \eqref{def:wcp}, \eqref{sdp}, and 
\eqref{Ritz}, we can write
\begin{eqnarray}
\label{ThetaAndRho}
{\rm i}{(\f,{\dot \theta}(t))}_{L^2}
-{(\nabla\f,\nabla\theta(t))}_{L^2}
&=&-{\rm i}{(\f,{\dot \rho}(t))}_{L^2}
+{(\f,v(\psi(t)-\psi_h(t)))}_{L^2}\nonumber\\
&&+  {(\f,f[\psi(t)]-f[\psi_h(t)])}_{L^2}.
\end{eqnarray}
Plugging $\f=\theta(t)\in S_h$ into \eqref{ThetaAndRho} and taking 
the imaginary part of the resulting equation, we get the differential 
inequality
\begin{eqnarray}
\label{DiffInequality1}
\frac12\,\frac{{\rm d}}{{\rm d}t}\,{\|\theta(t)\|}_{L^2}^2
&\le& \left({\|{\dot \rho}(t)\|}_{L^2}
+{\|v(\psi(t)-\psi_h(t))\|}_{L^2}
+{\|f[\psi(t)]-f[\psi_h(t)]\|}_{L^2}\right) 
{\|\theta(t)\|}_{L^2}\nonumber\\
&\le&  \left({\|{\dot \rho}(t)\|}_{L^2}
+c_1\left({\|\rho(t)\|}_{L^2}+{\|\theta(t)\|}_{L^2}\right) \right) 
{\|\theta(t)\|}_{L^2},
\end{eqnarray}
where we used the conservation laws \eqref{cp:conservation} and 
\eqref{sdp:conservation} and \eqref{cond_dd} to define the constant
\begin{eqnarray*}
c_1
:={\|v\|}_{L^\infty}
+2{\|\lambda\|}_{L^\infty}
{\|V\|}_{L^\infty(\R^d)}\left(\mM[\psi_0]+\mM[\psi_{0h}]\right).
\end{eqnarray*}
Using  $\epsilon>0$ to regularize the time derivative of 
${\|\theta\|}_{L^2}$ at $\theta=0$ by rewriting the l.h.s. of 
\eqref{DiffInequality1} as
$\frac12\frac{{\rm d}}{{\rm d}t}{\|\theta\|}_{L^2}^2
=\frac12\frac{{\rm d}}{{\rm d}t}({\|\theta\|}_{L^2}^2+\epsilon^2)$, we get
\begin{eqnarray}
\label{DiffInequality2}
\frac{{\rm d}}{{\rm d}t}\left({\|\theta\|}_{L^2}^2+\epsilon^2\right)^{1/2}
\le {\|{\dot \rho}(t)\|}_{L^2}
+c_1\left({\|\rho(t)\|}_{L^2}+{\|\theta(t)\|}_{L^2}\right),
\end{eqnarray}
where we used ${\|\theta(t)\|}_{L^2}\le \left({\|\theta\|}_{L^2}^2+\epsilon^2\right)^{1/2}$. Integrating \eqref{DiffInequality2} from 
$0$ to $t$, letting $\epsilon\to 0$, and applying Gr\"onwall's lemma to the resulting inequality, we find
\begin{eqnarray}
\label{Groenwall}
{\|\theta(t)\|}_{L^2}
&\le& {\|\theta(0)\|}_{L^2}
+\int_0^t\!{\rm d}s\,\big({\|{\dot \rho}(s)\|}_{L^2}
+c_1{\|\rho(s)\|}_{L^2}\big)\nonumber\\
&&+\,c_1\int_0^t\!{\rm d}s\,\Big(
{\|\theta(0)\|}_{L^2}
+\int_0^s\!{\rm d}u\,\big({\|{\dot \rho}(u)\|}_{L^2}
+c_1{\|\rho(u)\|}_{L^2}\big)\Big)\,
{\rm e}^{c_1(t-s)}.
\end{eqnarray}
In order to extract the factor $h^2$, 
we apply \eqref{RitzEstimate} and Assumption \ref{ass:initial} to get
\begin{eqnarray}
\label{rhoEstimate}
{\|\rho(t)\|}_{L^2}
&\le&C_R{\|\psi(t)\|}_{H^2} h^2,\\
\label{rhotEstimate}
{\|{\dot \rho}(t)\|}_{L^2}
&\le&C_R{\|{\dot \psi}(t)\|}_{H^2} h^2,\\
\label{theta0Estimate}
{\|\theta(0)\|}_{L^2}
&\le& \big(C_0+C_R{\|\psi_0\|}_{H^2}\big) h^2.
\end{eqnarray}
Plugging \eqref{rhoEstimate}, \eqref{rhotEstimate},
 and \eqref{theta0Estimate} into \eqref{Groenwall}, we finally arrive at
\begin{eqnarray*}
{\|\psi(t)-\psi_h(t)\|}_{L^2}
\le {\|\theta(t)\|}_{L^2}+{\|\rho(t)\|}_{L^2}
\le c_2(t) h^2,
\end{eqnarray*}
where the time dependent prefactor is defined by
\begin{eqnarray*}
c_2(t)
&:=&\big(C_0+C_R{\|\psi_0\|}_{H^2}\big) {\rm e}^{c_1 t} 
+C_R\int_0^t\!{\rm d}s \,\big(
{\|{\dot \psi}(s)\|}_{H^2}+ c_1{\|\psi(s)\|}_{H^2}\big)\\
&&+c_1 C_R \int_0^t\!{\rm d}s\int_0^s\!{\rm d}u\, 
\big({\|{\dot \psi}(u)\|}_{H^2}+c_1{\|\psi(u)\|}_{H^2}\big)
\,{\rm e}^{c_1(t-s)}
+C_R {\|\psi(t)\|}_{H^2}.
\end{eqnarray*}
Setting $C_E:=\max_{t\in[0,T]}c_2(t)$ brings the proof of Theorem 
\ref{thm:ErrorSD} to an end.\hfill$\Box$\vspace{5mm}

\begin{rem}
For the local case $f[\psi]=|\psi|^2\psi$, one replaces the original 
locally Lipschitz nonlinearity $f$ by a globally Lipschitz continuous 
nonlinearity
which coincides with $f$ in a given neighborhood of the solution $\psi$
of the continuum problem. One then first 
shows that the semidiscrete solution of the modified problem satisfies  
the desired $L^2$-error bound, and, second, that for $h$ sufficiently small,
the modified solution lies in the given neighborhood of $\psi$. But for
such $h$, the solution of the modified problem coincides with the solution
of the original problem, and, hence, the solution of the original problem
satisfies the desired $L^2$-error bound, too.
\end{rem}

\section{The fully discrete approximation}
\label{sec:fd}

In this section, we discretize the semidiscrete problem \eqref{sdp} in 
time. To this end, let us denote by $N\in\N$ 
the desired fineness of the time discretization with time discretization 
scale $\tau$ and its multiples $t_n$ for all $n=0,1,2,...,N$,
\begin{eqnarray}
\label{def:tau}
\tau:=\frac{T}{N},\quad 
t_n:=n\tau.
\end{eqnarray}
As mentioned in the Introduction, we will use two different
time discretization schemes of Crank-Nicholson type to approximate the
semidiscrete solution  $\psi_h$ of Theorem 
\ref{sdp:ExUniqueCons} at time $t_n$ by $\Psi_n\in\Psi$, where
\begin{eqnarray}
\label{def:Psi}
\Psi:=(\Psi_0,\Psi_1,...,\Psi_N)\in S_h^{\times(N{+}1)}.
\end{eqnarray}
These two schemes differ in the 
way of approximating the nonlinear term 
$g_V[|\psi|^2]\psi$ as follows. Let $\mN:=\{1,2,...,N\}$ and 
$\mN_0:=\mN\cup\{0\}$, and define
\begin{eqnarray}
\label{def:average}
\Psi_{n-1/2}
:=\frac{1}{2}\left(\Psi_{n}+\Psi_{n-1}\right),
\quad\forall n\in\mN.
\end{eqnarray}
Then, the first scheme implements the
one-step one-stage 
Gauss-Legendre 
Runge-Kutta method in which the nonlinear term is discretized by
\begin{eqnarray}
\label{gCoherent}
g_V[|\Psi_{n-1/2}|^2]\Psi_{n-1/2}.
\end{eqnarray}
In this method,  the mass $\mM[\Psi_n]$ is conserved under the discrete time
 evolution. The second scheme, introduced in \cite{DFP81} and applied in 
\cite{ADK91}, discretizes the nonlinear term by 
\begin{eqnarray}
\label{gIncoherent}
g_V[\tfrac12(|\Psi_{n}|^2+|\Psi_{n-1}|^2)]\Psi_{n-1/2}.
\end{eqnarray}
This method, in addition to the mass, also conserves the energy 
$\mH[\Psi_n]$ of the system. In the following, for convenience, 
we will call the first scheme {\it coherent} and the second one 
{\it incoherent}.

\subsection{Coherent scheme}
\label{sec:cfd}

In order to define what we mean by a coherent solution of the Hartree
initial-boundary value problem 
\eqref{def:cp}, we define
\begin{eqnarray*}
{\dot \Psi}_n
:=\frac{1}{\tau}\left(\Psi_{n}-\Psi_{n-1}\right), 
\quad\forall n\in\mN.
\end{eqnarray*}

\begin{definition}
\label{def:cfd}
Let Assumption \ref{ass:Vv} and \ref{ass:Sh} hold.
We call $\Psi\in S_h^{\times(N{+}1)}$ a coherent fully discrete solution of 
the Hartree
initial-boundary value problem \eqref{def:cp} with initial condition
$\psi_{0h}\in S_h$ if
\begin{eqnarray}
\label{cfd}
\left\{
\begin{array}{l}
{\rm i}{(\f,{\dot \Psi}_n)}_{L^2}
={(\nabla\f,\nabla\Psi_{n-1/2})}_{L^2}
+\,{(\f,v\Psi_{n-1/2})}_{L^2}
+{(\f,f[\Psi_{n-1/2}])}_{L^2},\\
\hspace{2.5cm}\forall\f\in S_h,\, \forall n\in\mN,\\
\Psi_0=\psi_{0h}.
\end{array}
\right.
\end{eqnarray}
\end{definition}

The coherent solution has the following conservation property.

\begin{proposition}
\label{prop:cfdMassConservation}
Let $\Psi\in S_h^{\times(N{+}1)}$ be  a coherent fully discrete solution of 
the Hartree
initial-boundary value problem \eqref{def:cp} in the sense of Definition 
\ref{def:cfd}. Then, the  mass of $\Psi$ is conserved under the 
discrete time evolution, 
\begin{eqnarray}
\label{cfdMassConservation}
\mM[\Psi_n]=\mM[\psi_{0h}],
\quad\forall n\in\mN_0.
\end{eqnarray}
\end{proposition}

{\it Proof}\, If we plug $\f=\Psi_{n-1/2}$ into \eqref{cfd} 
and take the imaginary part of the resulting equation,  we get
\begin{eqnarray*}
0=
{\rm Im}\,{\rm i}{(\Psi_{n-1/2},{\dot \Psi}_n)}_{L^2}
=\frac{1}{2\tau}\left(\mM[\Psi_{n}]-\mM[\Psi_{n-1}]\right).
\end{eqnarray*}
\hfill $\Box$

\begin{rem}
The energy $\mH[\Psi_n]$ of the coherent solution \eqref{cfd} is not
conserved under the discrete time evolution, see \cite{ADK91} and references 
therein, in particular \cite{Sanz-Serna} and \cite{Tourigny} for the local 
case with $d=1$.
\end{rem}

The question of existence and uniqueness of a coherent
solution is addressed in the following.

\begin{thm}
\label{prop:ExistenceCFD}
Let Assumptions \ref{ass:Vv}, \ref{ass:Sh}, and \ref{ass:inverse} hold,
and let the time discretization scale $\tau$ be sufficiently small.
Then, there exists a unique coherent fully discrete solution of the Hartree 
initial-boundary value problem \eqref{def:cp} in the sense of Definition 
\ref{def:cfd}.
\end{thm}

{\it Proof}\, Let $\phi\in S_h$ be given, and define the mapping 
$F_\phi: S_h\to S_h$ by
\begin{eqnarray}
\label{def:Fn}
{(\f,F_\phi[\psi])}_{L^2}
:={(\f,\phi)}_{L^2}
-\frac{{\rm i}\tau}{2}\, \big({(\nabla\f,\nabla\psi)}_{L^2}
+{(\f,v\psi)}_{L^2}
+{(\f,f[\psi])}_{L^2}\big),\quad
\forall \f\in S_h.
\end{eqnarray}
For some  $n\in\mN$, let the $n$-th component $\Psi_{n-1}$ of 
$\Psi\in S_h^{\times(N{+}1)}$ from \eqref{def:Psi} be given. Adding 
$2{\rm i}\,(\f,\Psi_{n-1})/\tau$ on 
both sides of \eqref{cfd}, we can rewrite \eqref{cfd} with the help of 
\eqref{def:Fn} in the form of a fixed point equation for $\Psi_{n-1/2}$,
\begin{eqnarray}
\label{cfdFixedPoint}
\Psi_{n-1/2}=F_{\Psi_{n-1}}[\Psi_{n-1/2}],
\end{eqnarray}
from which we retrieve the unknown component $\Psi_n$ by \eqref{def:average}. In order
 to construct the unique solution of \eqref{cfdFixedPoint}, we 
make use of Banach's fixed point theorem on the compact ball 
$\mB_{n-1}:=\{\psi\in S_h\,|\, {\|\psi\|}_{L^2}\le \mM[\Psi_{n-1}]^{1/2}+1\}$
in $S_h$. Using Assumption \ref{ass:inverse} and \eqref{cond_dd}, we get, 
for $\psi,\xi\in S_h$,
\begin{eqnarray}
\label{banach1}
|{(\f,F_{\Psi_{n-1}}[\psi]-F_{\Psi_{n-1}}[\xi])}_{L^2}|&&\nonumber\\
&&\hspace{-4.5cm}\le\frac\tau 2
\Big(C_B^2 h^{-2}
+{\|v\|}_{L^\infty} 
+2{\|\lambda\|}_{L^\infty}{\|V\|}_{L^\infty(\R^d)}
\big({\|\psi\|}_{L^2}^2+{\|\xi\|}_{L^2}^2\big)\Big)
{\|\psi-\xi\|}_{L^2} {\|\f\|}_{L^2}.
\end{eqnarray}
Plugging $\f=F_{\Psi_{n-1}}[\psi]-F_{\Psi_{n-1}}[\xi]$ into 
\eqref{banach1} and picking 
$\psi$ and $\xi$ from $\mB_{n-1}$, we find
\begin{eqnarray}
\label{banach2}
{\|F_{\Psi_{n-1}}[\psi]-F_{\Psi_{n-1}}[\xi]\|}_{L^2}
\le \frac{\alpha_{n-1}}{2}\,\tau {\|\psi-\xi\|}_{L^2},
\end{eqnarray}
where the constant $\alpha_{n-1}$ is defined, for all $n\in\mN$,  by
\begin{eqnarray}
\label{taun}
\alpha_{n-1}
:=
C_B^2 h^{-2}
+{\|v\|}_{L^\infty}
+4{\|\lambda\|}_{L^\infty}{\|V\|}_{L^\infty(\R^d)}(\mM[\Psi_{n-1}]^{1/2}+1)^2.
\end{eqnarray}
Let now $\alpha_{n-1}(\mM[\Psi_{n-1}]^{1/2}+1)\tau\le 1$. Then, it 
follows from \eqref{banach2} and \eqref{taun} that $F_{\Psi_{n-1}}$ maps 
$\mB_{n-1}$ into $\mB_{n-1}$ (set $\xi=0$ in \eqref{banach2}) and that 
$F_{\Psi_{n-1}}$ is a strict contraction on  $\mB_{n-1}$. Therefore, 
for such $\tau$, Banach's fixed point theorem implies the existence of 
a unique solution $\Psi_{n-1/2}\in\mB_{n-1}$ 
of the fixed point equation \eqref{cfdFixedPoint}. Moreover, due to the mass 
conservation \eqref{cfdMassConservation}, there exists 
no solution $\Psi_{n-1/2}$ of \eqref{cfdFixedPoint} with 
$\Psi_{n-1/2}\in S_h\setminus \mB_{n-1}$. Hence, the component 
$\Psi_n$ of the coherent solution 
exists and is unique for such $\tau$. Starting at $\Psi_0=\psi_{0h}$ and 
proceeding iteratively, we 
get all $n+1$ components of the coherent solution 
$\Psi\in S_h^{\times(N{+}1)}$. Moreover, again due to  
\eqref{cfdMassConservation}, we get a uniform bound on the size of the time 
discretization scale $\tau$, e.g. 
\begin{eqnarray*}
\alpha_0(\mM[\psi_{0h}]^{1/2}+1)\tau\le 1.
\end{eqnarray*}
\hfill$\Box$

\begin{rem}
Since $\alpha_0\ge C_B^2h^{-2}$, we have that $\tau\le C_B^{-2}h^2$,
where $C_B$ stems from  Assumption \ref{ass:inverse}.
\end{rem}


We next turn to the first of the two main assertions of the present paper 
which is the time quadratic accuracy estimate on the 
$L^2$-error of the coherent solution. 

\begin{thm}
\label{thm:cfderror}
Let Assumptions \ref{ass:Vv}, \ref{ass:Sh}, \ref{ass:AccuracyOrder}, and
\ref{ass:initial} hold, and let $\Psi\in 
S_h^{\times(N{+}1)}$ be the coherent solution from Theorem
\ref{prop:ExistenceCFD}. Then, 
there exists a constant $C_K>0$ s.t.
\begin{eqnarray*}
\max_{n\in\mN_0}{\|\psi(t_n)-\Psi_n\|}_{L^2}
\le 
C_K (\tau^2+h^2).
\end{eqnarray*}
\end{thm} 

\begin{rem}
The constant $C_K$ depends on higher Sobolev norms of the continuum
solution $\psi$.  These norms exist due to the regularity assertion in
Theorem \ref{thm:cs} (c).
\end{rem}

{\it Proof}\, Let $n\in\mN$ be fixed and define 
$\psi_n:=\psi(t_n)$ with $t_n$ from \eqref{def:tau}. As in the proof of 
Theorem 
\ref{thm:ErrorSD}, we decompose the difference to be estimated as
\begin{eqnarray}
\label{decomposition}
\psi_n-\Psi_n=\rho_n+\theta_n,
\end{eqnarray}
where $\rho_n$ and $\theta_n$ are again defined with the help of the 
Ritz projection from \eqref{Ritz} by
\begin{eqnarray}
\label{def:RhoTheta}
\rho_n:=(1-R_h)\psi_n,\quad 
\theta_n:=R_h\psi_n-\Psi_n.
\end{eqnarray}
Using Taylor's theorem in order to expand $\psi_n$ around $t=0$ up to zeroth 
order in $t_n$ and the estimate on the Ritz projection \eqref{RitzEstimate}, 
we immediately get 
\begin{eqnarray}
\label{RhoEstimate}
{\|\rho_n\|}_{L^2}
\le C_Rh^2\Big({\|\psi_0\|}_{H^2}
+\int_0^{t_n}\!\!{\rm d}t\,\,{\|{\dot \psi}(t)\|}_{H^2}\Big).
\end{eqnarray}
In order to estimate $\theta_n$, we want to extract suitable small 
differences from the expression
\begin{eqnarray}
\label{def:Ln}
L_{n,\f}:=
\frac{{\rm i}}{\tau}{(\f,\theta_{n}-\theta_{n-1})}_{L^2}
-\frac12 {(\nabla\f,\nabla (\theta_{n}+\theta_{n-1}))}_{L^2}
-\frac12 {(\f,v (\theta_{n}+\theta_{n-1}))}_{L^2}
\end{eqnarray}
which contains all the linear terms in \eqref{cfd} moved to 
the l.h.s. with $\Psi_n$ replaced by $\theta_n$. For this purpose,
we first plug the definition of $\theta_n$ into 
\eqref{def:Ln}, and then use the definition of the Ritz projection 
\eqref{Ritz} and the scheme \eqref{cfd} to get
\begin{eqnarray}
\label{Ln1}
L_{n,\f}&=&
\frac{{\rm i}}{\tau}\,{(\f,R_h(\psi_{n}-\psi_{n-1}))}_{L^2}
-\frac12 {(\nabla\f,\nabla (\psi_{n}+\psi_{n-1}))}_{L^2}
-\frac12 {(\f,v R_h(\psi_{n}+\psi_{n-1}))}_{L^2}\nonumber\\
&&-\,{(\f,f[\Psi_{n-1/2}])}_{L^2}.
\end{eqnarray}
Rewriting the first term on the r.h.s. of \eqref{Ln1} with the help of
the continuum solution satisfying the weak formulation \eqref{def:wcp}, 
we have 
\begin{eqnarray}
\label{Rh}
\frac{{\rm i}}{\tau}\,{(\f,R_h(\psi_{n}-\psi_{n-1}))}_{L^2}
&=&\frac{{\rm i}}{\tau}\,{(\f,(R_h-1)(\psi_{n}-\psi_{n-1}))}_{L^2}
+{\rm i}\,{(\f,\tfrac{1}{\tau}(\psi_{n}-\psi_{n-1})-{\dot \psi}_{n-1/2})}_{L^2}
\nonumber\\
&&+\,{(\nabla\f,\nabla\psi_{n-1/2})}_{L^2}
+{(\f,v\psi_{n-1/2})}_{L^2}
+{(\f,f[\psi_{n-1/2}])}_{L^2},
\end{eqnarray}
where we used the notations $\psi_{n-1/2}:=\psi(t_n-\tau/2)$ and
 ${\dot \psi}_{n-1/2}:={\dot \psi}(t_n-\tau/2)$. Plugging \eqref{Rh} into
\eqref{Ln1}, we can express  $L_{n,\f}$ in the form
\begin{eqnarray}
\label{def:omegas}
L_{n,\f}=\sum_{j=1}^6\,{(\f,\omega_n^{(j)})}_{L^2},
\end{eqnarray}
where the functions $\omega_n^{(j)}$ with $j=1,...,6$ are defined by
\begin{eqnarray*}
\omega_n^{(1)}
&:=&\frac{{\rm i}}{\tau}(R_h-1)(\psi_{n}-\psi_{n-1}),\\
\omega_n^{(2)}
&:=&{\rm i}\big(\tfrac{1}{\tau}(\psi_{n}-\psi_{n-1})
-{\dot \psi}_{n-1/2}\big),\\
\omega_n^{(3)}
&:=&\Delta\big(\tfrac12\left(\psi_{n}+\psi_{n-1}\big)
-\psi_{n-1/2}\right),\\
\omega_n^{(4)}
&:=&v\left(\psi_{n-1/2}-\tfrac12\left(\psi_{n}+\psi_{n-1}\right)\right),\\
\omega_n^{(5)}
&:=&\frac12 \,v\left(1-R_h\right)\left(\psi_{n}+\psi_{n-1}\right),\\
\omega_n^{(6)}
&=& f[\psi_{n-1/2}]-f[\Psi_{n-1/2}],
\end{eqnarray*}
and, for $\omega_n^{(3)}$, we used again 
${(\nabla\f,\nabla\psi)}_{L^2}={(\f,-\Delta\psi)}_{L^2}$ for all 
$\psi\in H^2\cap H^1_0$. Plugging $\f=(\theta_{n}+\theta_{n-1})/2$ into 
\eqref{def:Ln} and \eqref{def:omegas},  and
taking the imaginary part of the resulting equation, we 
get\footnote{If $\theta_{n}=0$, we are left with \eqref{RhoEstimate}.}
\begin{eqnarray}
\label{thetaEstimate1}
{\|\theta_{n}\|}_{L^2}
\le {\|\theta_{n-1}\|}_{L^2}
+\tau\sum_{j=1}^6 {\|\omega_n^{(j)}\|}_{L^2}.
\end{eqnarray}
Let us next estimate the terms ${\|\omega_n^{(j)}\|}_{L^2}$ for all 
$j=1,...,6$. 
For $\omega_n^{(1)}$, we 
expand $\psi_{n}$ around $t=t_{n-1}$ up to zeroth order in $\tau$ and 
use \eqref{RitzEstimate} s.t.
\begin{eqnarray}
\label{omega1Estimate}
{\|\omega_n^{(1)}\|}_{L^2}
\le {C_Rh^2}\tau^{-1}\int_{t_{n-1}}^{t_{n}}\!\!{\rm d}t\,\,
{\|{\dot \psi}(t)\|}_{H^2}.
\end{eqnarray}
For $\omega_n^{(2)}$,  we expand $\psi_{n-1}$ and $\psi_{n}$ around 
$t=t_n-\tau/2$ up to second order in $\tau/2$,
\begin{eqnarray}
\label{omega2Estimate}
{\|\omega_n^{(2)}\|}_{L^2}
&\le& \frac{1}{2\tau}\Big(\int_{t_{n-1}}^{t_n-\tau/2}\!\!{\rm d}t\,\,
(t_{n-1}-t)^2{\|{\dddot \psi}(t)\|}_{L^2}
+\int_{t_n-\tau/2}^{t_{n}}\!\!{\rm d}t\,\,
(t_{n}-t)^2{\|{\dddot \psi}(t)\|}_{L^2}\Big)\nonumber\\
&\le& \frac{\tau}{8}\int_{t_{n-1}}^{t_{n}}\!\!{\rm d}t\,\,
{\|{\dddot \psi}(t)\|}_{L^2}.
\end{eqnarray}
Analogously, for $\omega_n^{(3)}$ and $\omega_n^{(4)}$, we expand $\psi_{n-1}$ 
and $\psi_{n}$ around $t=t_n-\tau/2$ up to first order in $\tau/2$,
\begin{eqnarray}
\label{omega3Estimate}
{\|\omega_n^{(3)}\|}_{L^2}
&\le& \frac{\tau}{4}\int_{t_{n-1}}^{t_{n}}\!\!{\rm d}t\,\,
{\|\Delta{\ddot \psi}(t)\|}_{L^2},\\
\label{omega4Estimate}
{\|\omega_n^{(4)}\|}_{L^2}
&\le& \frac{\tau}{4}\,{\|v\|}_{L^\infty}\int_{t_{n-1}}^{t_{n}}\!\!{\rm d}t\,\,
{\|{\ddot \psi}(t)\|}_{L^2}.
\end{eqnarray}
For 
$\omega_n^{(5)}$, expanding $\psi_{n-1}$ and $\psi_{n}$ around $t=0$ up to zeroth order in time, we get,  analogously to the estimate of $\omega_n^{(1)}$, 
\begin{eqnarray}
\label{omega5Estimate}
{\|\omega_n^{(5)}\|}_{L^2}
\le C_R h^2 {\|v\|}_{L^\infty} \Big({\|\psi_0\|}_{H^2}
+\int_{0}^{t_{n}}\!\!{\rm d}t\,\,{\|{\dot \psi}(t)\|}_{H^2}\Big).
\end{eqnarray}
Finally, for $\omega_n^{(6)}$, we apply the local Lipschitz continuity 
\eqref{cond_dd} to get
\begin{eqnarray}
\label{omega6Estimate1}
{\|\omega_n^{(6)}\|}_{L^2}
\le c_1\,{\|\psi_{n-1/2}-\Psi_{n-1/2}\|}_{L^2},
\end{eqnarray}
where we used the continuum mass conservation \eqref{cp:conservation} 
and the coherent fully discrete mass 
conservation \eqref{cfdMassConservation} to define the constant
\begin{eqnarray}
\label{def:c1}
c_1
:=2 {\|\lambda\|}_{L^\infty}{\|V\|}_{L^\infty(\R^d)}
\left(\mM[\psi_0]+\mM[\psi_{0h}]\right).
\end{eqnarray}
Since we want to reinsert the decomposition \eqref{decomposition}
 into the r.h.s. of \eqref{omega6Estimate1}, we write
\begin{eqnarray}
\label{insert}
{\|\psi_{n-1/2}-\Psi_{n-1/2}\|}_{L^2}
&\le& 
{\|\psi_{n-1/2}-\tfrac12(\psi_{n}+\psi_{n-1})\|}_{L^2}\nonumber\\
&&
+\,\frac12 \left({\|\rho_{n-1}\|}_{L^2}+{\|\rho_{n}\|}_{L^2}
+{\|\theta_{n-1}\|}_{L^2}+{\|\theta_{n}\|}_{L^2}\right).
\end{eqnarray}
Plugging the estimates 
\eqref{RhoEstimate}, \eqref{omega1Estimate} to
\eqref{omega6Estimate1}, and \eqref{insert} 
into \eqref{thetaEstimate1}, we find
\begin{eqnarray}
\label{ThetaClosed}
{\|\theta_{n}\|}_{L^2}
{\le \|\theta_{n-1}\|}_{L^2}
+\big(A^{(1)}_n+\tau A^{(2)}_n\big)h^2
+A^{(3)}_n\tau^2
+\frac{c_1}{2}\,\tau \left({\|\theta_{n-1}\|}_{L^2}
+{\|\theta_{n}\|}_{L^2}\right),
\end{eqnarray}
where the first term on the r.h.s. of \eqref{insert} was estimated as in 
$\omega_n^{(3)}$ or $\omega_n^{(4)}$, and
\begin{eqnarray}
\label{def:A1}
A^{(1)}_n
&:=&C_R\int_{t_{n-1}}^{t_{n}}\!\! {\rm d}t\,\,{\|\dot \psi(t)\|}_{H^2},\\
\label{def:A2}
A^{(2)}_n
&:=&C_R \left(c_1+{\|v\|}_{L^\infty}\right)
\Big({\|\psi_0\|}_{H^2}
+\int_{0}^{t_{n}}\!\!{\rm d}t\,\,{\|\dot\psi(t)\|}_{H^2}\Big),\\
\label{def:A3}
A^{(3)}_n
&:=&\frac14\left(c_1+{\|v\|}_{L^\infty}\right)
\int_{t_{n-1}}^{t_{n}}\!\!{\rm d}t\,\,
{\|\ddot\psi(t)\|}_{L^2}
+\frac18\int_{t_{n-1}}^{t_{n}}\!\!{\rm d}t\,\,
{\|\dddot\psi(t)\|}_{L^2}
+\,\frac14\int_{t_{n-1}}^{t_{n}}\!\!{\rm d}t\,\,
{\|\Delta\ddot\psi(t)\|}_{L^2}.\quad
\end{eqnarray}
If we choose the time discretization scale $\tau$ to be small enough, 
e.g. $c_1\tau\le 1$,
we can construct the following recursive bound on 
${\|\theta_{n}\|}_{L^2}$ from inequality \eqref{ThetaClosed},
\begin{eqnarray}
\label{theta_n}
{\|\theta_{n}\|}_{L^2}
\le B^{(1)}{\|\theta_{n-1}\|}_{L^2}+B^{(2)}_n,
\end{eqnarray}
where we used that $1/(1-c_1\tau/2)\le 1+c_1\tau$ if $c_1\tau\le 1$ to define
\begin{eqnarray*}
B^{(1)}
&:=&(1+c_1\tau)^2,\\
B^{(2)}_n
&:=&(1+c_1\tau)\big(\!\big(A^{(1)}_n
+\tau A^{(2)}_n\big)h^2+A^{(3)}_n\tau^2\big).
\end{eqnarray*}
Therefore, if we iterate the bound \eqref{theta_n} 
until we arrive at ${\|\theta_0\|}_{L^2}$, we get 
\begin{eqnarray}
\label{theta_n2}
{\|\theta_n\|}_{L^2}
&\le& \left(B^{(1)}\right)^n{\|\theta_0\|}_{L^2}
+\sum_{k=1}^{n}\left(B^{(1)}\right)^{n-k} B^{(2)}_ {k}\nonumber\\
&\le& c_2\Big({\|\theta_0\|}_{L^2}
+\sum_{k=1}^{n}\big(\big(A^{(1)}_{k}
+\tau A^{(2)}_{k}\big)h^2
+A^{(3)}_{k}\tau^2 \big)\Big),
\end{eqnarray}
where, on the second line of \eqref{theta_n2}, we first extract the global factor  $(B^{(1)})^n$ which can then be estimated as 
$(B^{(1)})^n\le (1+c_1 T/N)^{2N}\le c_2$ with the definition
\begin{eqnarray}
\label{def:c2}
c_2
:={\rm e}^{2c_1 T}.
\end{eqnarray}
It remains to estimate ${\|\theta_0\|}_{L^2}$ on the r.h.s. of 
\eqref{theta_n2}.
This is again done by using the estimate on the Ritz projection 
\eqref{RitzEstimate},
\begin{eqnarray}
\label{theta_0}
{\|\theta_0\|}_{L^2}
&\le& {\|\psi_0-\psi_{0h}\|}_{L^2}
+{\|(R_h-1)\psi_0\|}_{L^2}\nonumber\\
&\le& {\|\psi_0-\psi_{0h}\|}_{L^2}
+C_Rh^2{\|\psi_0\|}_{H^2}.
\end{eqnarray}
Hence, with estimate \eqref{RhoEstimate} on ${\|\rho_n\|}_{L^2}$ and 
the estimates  \eqref{theta_n2} and \eqref{theta_0} on ${\|\theta_n\|}_{L^2}$ 
in the decomposition \eqref{def:RhoTheta}, 
taking the maximum over all times, we finally arrive at
\begin{eqnarray}
\label{cfd:max}
\max_{n\in\mN_0}{\|\psi_n-\Psi_n\|}_{L^2}
\le c_2{\|\psi_0-\psi_{0h}\|}_{L^2}+c_3 h^2+c_4\tau^2,
\end{eqnarray}
where the constants $c_3$ and $c_4$ are defined by
\begin{eqnarray}
\label{cfd:c3}
c_3
&:=&C_R\left(1+c_2\left(1+T 
\left(c_1
+{\|v\|}_{L^\infty}\right)\right)\right)
\Big({\|\psi_0\|}_{H^2}+\int_0^T\!\!{\rm d}t\,\,{\|\dot\psi(t)\|}_{H^2}\Big),
\\
\label{cfd:c4}
c_4
&:=&\frac{c_2}{4}\,\Big((c_1+{\|v\|}_{L^\infty})
\int_0^T\!\!{\rm d}t\,\,{\|\ddot\psi(t)\|}_{L^2}
+\frac12 \int_0^T\!\!{\rm d}t\,\,{\|\dddot\psi(t)\|}_{L^2}
+\int_0^T\!\!{\rm d}t\,\,{\|\Delta\ddot\psi(t)\|}_{L^2}\Big).
\end{eqnarray}
The constants $c_1$ and $c_2$ are given in \eqref{def:c1} and \eqref{def:c2}, 
respectively. Using Assumption \ref{ass:initial} and setting 
$C_K:=\max\{c_2C_0+c_3,c_4\}$ brings the proof of Theorem \ref{thm:cfderror}
 to an end.\hfill $\Box$ \vspace{5mm}

\begin{rem}
\label{rem:bound}
We can compute an explicit bound on the integrands in \eqref{cfd:c3} and 
\eqref{cfd:c4} on any finite time interval. As an example, for the last term 
in the constant $c_4$ from
\eqref{cfd:c4}, we have\footnote{See Theorem \ref{thm:cs} (c) 
and \cite[p.299]{RS2}.}
\begin{eqnarray*}
{\|\Delta\ddot\psi(t)\|}_{L^2}
\le {\|\Delta^3\psi(t)\|}_{L^2}
+{\|\Delta^2J[\psi(t)]\|}_{L^2}
+{\|\Delta \tfrac{{\rm d}}{{\rm d}t}J[\psi(t)]\|}_{L^2}.
\end{eqnarray*}
The first term is exponentially bounded in time using the conditions
from Theorem \ref{thm:cs} (a) and (b) 
and  Gr\"onwall's
lemma on the Duhamel 
integral form  of the differential equation \eqref{cp},
\begin{eqnarray*}
{\|\Delta^3\psi(t)\|}_{L^2}
\le {\|\Delta^3\psi_0\|}_{L^2}\, {\rm e}^{C(\mM[\psi_0])t},
\end{eqnarray*}
where the constant $C$ stems from \eqref{cond_cc}.\footnote{See also 
\cite[p.300]{RS2}. In contradistinction to the general case from
Theorem \ref{thm:cs} (a), the growth rate $C$ from \eqref{cond_cc} 
only depends on the mass of the initial condition (and on $v$, $V$, 
and $\lambda$, of course).}
The second term is again bounded 
due to \eqref{cond_c}. Finally, the third 
term is bounded due to equation \eqref{deriv} for the time derivative of the 
nonlinear term $J[\psi(t)]$ and the corresponding estimates \eqref{est1}
and \eqref{est2}.
\end{rem}

\subsection{Incoherent scheme}

As described at the beginning of Section \ref{sec:fd}, we also study a second 
discretization scheme which approximizes the nonlinear term 
$g_V[|\psi|^2]\psi$ not by \eqref{gCoherent} but rather 
by the expression \eqref{gIncoherent}.

\begin{definition}
\label{def:ifd}
Let Assumptions \ref{ass:Vv} and \ref{ass:Sh} hold.
We call $\Psi\in S_h^{\times (N+1)}$ an incoherent fully discrete solution of 
the Hartree initial-boundary value problem \eqref{def:cp} with initial
condition $\psi_{0h}\in S_h$ if
\begin{eqnarray}
\label{ifd}
\left\{
\begin{array}{l}
{\rm i}\,{(\f,\dot\Psi_n)}_{L^2}
={(\nabla\f,\nabla\Psi_{n-1/2})}_{L^2}
+{(\f,v\Psi_{n-1/2})}_{L^2}+
{(\f,\tfrac12\lambda g_V[|\Psi_n|^2+|\Psi_{n-1}|^2]\Psi_{n-1/2})}_{L^2},\\
\hspace{2.5cm}\forall\f\in S_h,\,\,\forall n\in\mN,\\
\Psi_0=\psi_{0h}. 
\end{array}
\right.
\end{eqnarray}
\end{definition}

The incoherent solution has the following conservation properties.

\begin{proposition}
\label{prop:ifdMassConservation}
Let $\Psi\in S_h^{\times (N+1)}$ be an incoherent fully discrete solution of 
the Hartree initial-boundary value problem \eqref{def:cp} in the sense of 
Definition \ref{def:ifd}, and let Assumption \ref{ass:even} hold.
Then, the mass and the energy of $\Psi$  are conserved under the discrete
time evolution, 
\begin{eqnarray}
\label{ifd:conservation}
\mM[\Psi_n]&=&\mM[\psi_{0h}],\quad
\forall n\in\mN_0,\\
\mH[\Psi_n]&=&\mH[\psi_{0h}],\quad
\forall n\in\mN_0.\nonumber
\end{eqnarray}
\end{proposition}

{\it Proof}\, Plugging $\f=\Psi_{n-1/2}$ into \eqref{ifd} and taking the
 imaginary part of the resulting equation leads to the mass conservation
as in the proof of Proposition \ref{prop:cfdMassConservation}. 
In order to prove the energy conservation, we plug $\f=\dot\Psi_n$ into 
\eqref{ifd} and take the real part of the resulting equation. Using 
that ${(\Psi_n,\lambda g_V[|\Psi_{n-1}|^2]\Psi_n)}_{L^2}
={(\Psi_{n-1},\lambda g_V[|\Psi_n|^2]\Psi_{n-1})}_{L^2}$ due to Assumption 
\ref{ass:even}, we get
\begin{eqnarray*}
0
={\rm Re}\,\,{\rm i}{(\dot \Psi_n,\dot \Psi_n)}_{L^2}
=\frac{1}{2\tau}\left(\mH[\Psi_n]-\mH[\Psi_{n-1}]\right).
\end{eqnarray*}
\hfill $\Box$ \vspace{5mm}

We next turn to the proof of existence-uniqueness of the incoherent 
solution.

\begin{thm}
\label{thm:ifdEx}
Let Assumptions \ref{ass:Vv}, \ref{ass:Sh}, and \ref{ass:inverse} hold,
and let the time discretization scale $\tau$ be sufficiently small.
Then, there exists a unique incoherent fully discrete solution of the 
Hartree initial-
boundary value problem \eqref{def:cp} in the sense of Definition 
\ref{def:ifd}.
\end{thm}

{\it Proof}\, The proof for the incoherent solution is 
analogous to the proof of the coherent solution. 
Let $\phi\in S_h$ be given, and define the mapping $G_\phi: S_h\to S_h$ by
\begin{eqnarray}
\label{def:Gphi}
{(\f,G_\phi[\psi])}_{L^2}
&:=&{(\f,\phi)}_{L^2}\nonumber\\
&&\hspace{-2.5cm}-\frac{{\rm i}\tau}{2}\left( 
{(\nabla\f,\nabla\psi)}_{L^2}
+{(\f,v\psi)}_{L^2}
+{(\f,\tfrac12\lambda g_V[|2\psi-\phi|^2+|\phi|^2]\psi)}_{L^2}
\right),\quad
\forall\f\in S_h.
\end{eqnarray}
For some $n\in\mN$, let the $n$-th component $\Psi_{n-1}$ of 
$\Psi\in S_h^{\times (N+1)}$ from \eqref{def:Psi} be given.
Adding $2{\rm i}(\f,\Psi_{n-1})/\tau$ on both sides of \eqref{ifd},
we rewrite \eqref{ifd} with the help of \eqref{def:Gphi} in the form
of a fixed point equation for $\Psi_{n-1/2}$,
\begin{eqnarray*}
\Psi_{n-1/2}
=G_{\Psi_{n-1}}[\Psi_{n-1/2}].
\end{eqnarray*}
In order to make 
use of Banach's fixed point theorem as in the proof of 
Theorem \ref{prop:ExistenceCFD}, we show that $G_{\Psi_{n-1}}$ 
maps the compact
Ball $\mB_{n-1}:=\{\psi\in S_h\,|\, {\|\psi\|}_{L^2}\le 
\mM[\Psi_{n-1}]^{1/2}+1\}$ 
into itself and that $G_{\Psi_{n-1}}$ is a strict contraction on 
$\mB_{n-1}$. To this end, we write
\begin{eqnarray}
\label{DiffGn1}
|{(\f,G_{\Psi_{n-1}}[\psi]-G_{\Psi_{n-1}}[\xi])}_{L^2}|\nonumber\\
&&\hspace{-4cm}\le\frac{\tau}{2}\Big(
{\|\nabla(\psi-\xi)\|}_{L^2} {\|\nabla\f\|}_{L^2}
+{\|v\|}_{L^\infty} {\|\psi-\xi\|}_{L^2} {\|\f\|}_{L^2}
+\,\frac12\,A{\|\f\|}_{L^2}
\Big),
\end{eqnarray}
where,  with the help of \eqref{est}, \eqref{dc1}, and 
$||z|^2-|w|^2|\le |z+w||z-w|$ for all $z,w\in\C$,
the third term $A$ on the r.h.s. of \eqref{DiffGn1} can be estimated as
\begin{eqnarray*}
A
&:=&{\|\lambda g_V[|2\psi-\Psi_{n-1}|^2+|\Psi_{n-1}|^2]\psi
-\lambda g_V[|2\xi-\Psi_{n-1}|^2+|\Psi_{n-1}|^2]\xi\|}_{L^2}\nonumber\\
&\le& 8{\|\lambda\|}_{L^\infty}{\|V\|}_{L^\infty(\R^d)}\,
\big({\|\psi\|}_{L^2}^2
+{\|\xi\|}_{L^2}^2
+{\|\Psi_{n-1}\|}_{L^2}^2\big)
 {\|\psi-\xi\|}_{L^2}.
\end{eqnarray*}
Hence, plugging $\f=G_{\Psi_{n-1}}[\psi]-G_{\Psi_{n-1}}[\xi]$ into 
\eqref{DiffGn1}, we get for $\psi,\xi\in\mB_{n-1}$ using Assumption 
\ref{ass:inverse},
\begin{eqnarray*}
{\|G_{\Psi_{n-1}}[\psi]-G_{\Psi_{n-1}}[\xi]\|}_{L^2}
\le \frac{\alpha_{n-1}}{2}\,\tau\,{\|\psi-\xi\|}_{L^2},
\end{eqnarray*}
where $\alpha_{n-1}
:=C_B^2h^{-2}
+{\|v\|}_{L^\infty} 
+12{\|\lambda\|}_{L^\infty}{\|V\|}_{L^\infty(\R^d)}\,(\mM[\Psi_{n-1}]^{1/2}+1)^2$ like in the coherent
scheme \eqref{taun}.
Therefore, we arrive at the claim as in the proof of Proposition
\ref{prop:ExistenceCFD} using the mass conservation from 
\eqref{ifd:conservation}, i.e. the incoherent solution exists and is 
unique if the time discretization scale  $\tau$ is sufficiently small, e.g.
$\alpha_0(\mM[\psi_{0h}]^{1/2}+1)\tau\le 1$. 
\hfill $\Box$ \vspace{5mm}

Finally, we also provide a time quadratic accuracy estimate on the 
$L^2$-error of the incoherent solution. Again, the proof is analogous
to the corresponding proof for the coherent solution from Theorem
\ref{thm:cfderror}.

\begin{thm}
\label{thm:ifderror}
Let Assumptions \ref{ass:Vv}, \ref{ass:Sh}, \ref{ass:AccuracyOrder}, and  
\ref{ass:initial} hold, and let $\Psi\in 
S_h^{\times(N{+}1)}$ be the incoherent solution from Theorem 
\ref{thm:ifdEx}. Then, 
there exists a constant $C_I>0$ s.t.
\begin{eqnarray*}
\max_{n\in \mN_0}{\|\psi(t_n)-\Psi_n\|}_{L^2}
\le 
C_I (\tau^2+h^2).
\end{eqnarray*}
\end{thm}

{\it Proof}\, As in \eqref{decomposition} and \eqref{def:RhoTheta}, 
we make use of the decomposition  $\psi_n-\Psi_n=\rho_n+\theta_n$,
and  we estimate $\rho_n$ again by \eqref{RhoEstimate}. In order to estimate 
${\|\theta_n\|}_{L^2}$, we define $L_{n,\f}$
as in \eqref{def:Ln}. Then, everything in equation \eqref{Ln1} remains
unchanged up to the last term which 
is replaced by the expression
$-(\f,\tfrac12 \lambda g_V[|\Psi_{n}|^2+|\Psi_{n-1}|^2]\Psi_{n-1/2})$. 
Using again 
\eqref{Rh}, we  can rewrite $L_{n,\f}$ as in equation \eqref{def:omegas},
where the terms $\omega_n^{(j)}$ remain unchanged for all $j=1,...,5$,
whereas the term $\omega_n^{(6)}$ now has the form
\begin{eqnarray*}
\omega_n^{(6)}
&:=&f[\psi_{n-1/2}]-\lambda g_V[\tfrac12(|\Psi_n|^2+|\Psi_{n-1}|^2)]
\Psi_{n-1/2}\\
&=&\omega_{n,1}^{(6)}
+\omega_{n,2}^{(6)}
+\omega_{n,3}^{(6)},
\end{eqnarray*}
where we use the same notation as introduced after \eqref{Rh} to define
\begin{eqnarray*}
\omega_{n,1}^{(6)}
&:=&\lambda g_V[|\psi_{n-1/2}|^2](\psi_{n-1/2}-\tfrac12 (\psi_n+\psi_{n-1})),\\
\omega_{n,2}^{(6)}
&:=&\frac12\, \lambda g_V[|\psi_{n-1/2}|^2-\tfrac12(|\psi_n|^2+|\psi_{n-1}|^2)]
(\psi_n+\psi_{n-1}),\\
\omega_{n,3}^{(6)}
&:=&\frac12\, \lambda g_V[\tfrac12(|\psi_n|^2+|\psi_{n-1}|^2)]
(\psi_n+\psi_{n-1})
- \lambda g_V[\tfrac12(|\Psi_n|^2+|\Psi_{n-1}|^2)]\Psi_{n-1/2}.
\end{eqnarray*}
For $\omega_{n,1}^{(6)}$, expanding $\psi_{n-1}$ and $\psi_{n}$ around 
$t=t_n-\tau/2$ up to first
order in $\tau/2$, we get , using \eqref{est} and the mass 
conservation \eqref{cp:conservation},
\begin{eqnarray}
\label{omega6_1}
{\|\omega_{n,1}^{(6)}\|}_{L^2}
\le \frac{\tau}{4}\,{\|\lambda\|}_{L^\infty}{\|V\|}_{L^\infty(\R^d)} 
\mM[\psi_0] 
\int_{t_{n-1}}^{t_n}\!{\rm d}t\,\, {\|\ddot\psi(t)\|}_{L^2}.
\end{eqnarray}
In order to estimate $\omega_{n,2}^{(6)}$, we expand
$\psi_{n-1}$ and $\psi_{n}$ around $t=t_n-\tau/2$ up to first
order in $\tau/2$ and get similarly
\begin{eqnarray}
\label{omega6_2}
{\|\omega_{n,2}^{(6)}\|}_{L^2}
&\le& {\|\lambda\|}_{L^\infty}{\|V\|}_{L^\infty(\R^d)  } \mM[\psi_0]^{1/2}{\| |\psi_{n-1/2}|^2
-\tfrac12 (|\psi_n|^2+|\psi_{n-1}|^2)\|}_{L^1}\nonumber\\
&\le& a^{(0)}(a_n^{(1)} \tau+a_n^{(2)} \tau^2+a_n^{(3)} \tau^3),
\end{eqnarray}
where we define
\begin{eqnarray*}
a^{(0)}
&:=& {\|\lambda\|}_{L^\infty}{\|V\|}_{L^\infty(\R^d)    } \mM[\psi_0]^{1/2},\\
a_n^{(1)}
&:=& \mM[\psi_0]\int_{t_{n-1}}^{t_n}\!{\rm d}t\,\, 
{\|\ddot\psi(t)\|}_{L^2}\\
a_n^{(2)}
&:=&\frac12 \,{\|\dot\psi_{n-1/2}\|}_{L^2}\Big({\|\dot\psi_{n-1/2}\|}_{L^2}
+\int_{t_{n-1}}^{t_n}\!{\rm d}t\,\, 
{\|\ddot\psi(t)\|}_{L^2}\Big)\\
a_n^{(3)}
&:=& \frac{1}{8}\int_{t_{n-1}}^{t_n}\!{\rm d}t\,\, 
{\|\ddot\psi(t)\|}_{L^2}^2.
\end{eqnarray*}
For $\omega_{n,3}^{(6)}$, using in particular again 
$||z|^2-|w|^2|\le |z+w| |z-w|$ for all 
$z,w\in\C$ and the decomposition
\eqref{decomposition}, we get 
\begin{eqnarray}
\label{omega6_3}
{\|\omega_{n,3}^{(6)}\|}_{L^2}
&\le& \frac14\, {\|\lambda g_V[|\psi_n|^2+|\psi_{n-1}|^2](\psi_n-\Psi_n+\psi_{n-1}-
\Psi_{n-1})\|}_{L^2}\nonumber\\
&&+\frac14\,{\|\lambda g_V[|\psi_n|^2-|\Psi_n|^2+|\psi_{n-1}|^2-|\Psi_{n-1}|^2]
(\Psi_n+\Psi_{n-1})\|}_{L^2}\nonumber\\
&\le&  {\|\lambda\|}_{L^\infty} {\|V\|}_{L^\infty(\R^d)}(\mM[\psi_0]+\mM[\psi_{0h}])
(\|\rho_n\|+\|\rho_{n-1}\|+\|\theta_n\|+\|\theta_{n-1}\|).
\end{eqnarray}
Therefore, plugging the estimates \eqref{omega6_1}, \eqref{omega6_2}, and 
\eqref{omega6_3} into \eqref{thetaEstimate1}, we again find the closed 
inequality
\eqref{ThetaClosed}, the coefficients $A_n^{(1)}$ and $A_n^{(2)}$ having the 
same form as in \eqref{def:A1} and \eqref{def:A2}, respectively.
Using estimate \eqref{omega6_1}  on $\omega_{n,1}^{(6)}$, we see
that the coefficient $A_n^{(3)}$ in the incoherent case contains all the 
terms from \eqref{def:A3} of the coherent case with $c_1$ replaced by
$a^{(0)}$, plus an additional 
term of the form $a^{(0)}(a_n^{(1)}+a_n^{(2)}\tau+a_n^{(3)}\tau^2)$ 
which is due to the estimate \eqref{omega6_2} of $\omega_{n,2}^{(6)}$.
 Plugging the coefficients $A_n^{(1)}$,
 $A_n^{(2)}$, and $A_n^{(3)}$ into the iterated bound \eqref{theta_n2} and 
using estimate \eqref{theta_0} on $\theta_0$, we again get an estimate of
 the form \eqref{cfd:max} where the constant $c_3$ has the same form as 
in \eqref{cfd:c3} whereas the constant $c_4$, compared to \eqref{cfd:c4},
 now looks like 
\begin{eqnarray}
\label{ifd:c4}
c_4
&:=&
\frac{c_2}{4} \Big((a^{(0)}+{\|v\|}_{L^\infty})
\int_0^T\!\!{\rm d}t\,\,{\|\ddot\psi(t)\|}_{L^2}
+\frac12 \int_0^T\!\!{\rm d}t\,\,{\|\dddot \psi(t)\|}_{L^2}
+\int_0^T\!\!{\rm d}t\,\,{\|\Delta\ddot\psi(t)\|}_{L^2}\Big)\nonumber\\
&&
+c_2a^{(0)}\big(\mM[\psi_0]
+\tfrac{1}{2}\,T\max_{t\in[0,T]}{\|\dot\psi(t)\|}_{L^2}\big)\int_0^T\!\!{\rm d}t\,\,{\|\ddot\psi(t)\|}_{L^2}\nonumber\\
&&
+\frac{c_2a^{(0)}}{2}T\max_{t\in[0,T]}{\|\dot\psi(t)\|}_{L^2}^2
+\frac{c_2a^{(0)}}{8}T^2 \int_0^T\!\!{\rm d}t\,\,{\|\ddot\psi(t)\|}_{L^2}^2.
\end{eqnarray}
Herewith, as in the proof of Theorem \ref{thm:cfderror}, we arrive
at the assertion. \hfill $\Box$ \vspace{5mm}

\begin{rem}
Using estimates as in Remark \ref{rem:bound}, we can again bound the 
constants in
\eqref{ifd:c4} explicitly.
\end{rem}


\end{document}